\documentclass[12pt,letterpaper]{elsarticle}
\usepackage[utf8]{inputenc}
\usepackage{amsmath}
\usepackage{amssymb}
 \usepackage{graphicx}
\usepackage{color}
\usepackage{url}
\usepackage{array}
\usepackage{caption} 
 \usepackage{setspace}
 \usepackage{lineno}
 
\newtheorem{definition}{Definition}
\newtheorem{lemma}{Lemma}
\newtheorem{proposition}{Proposition}
\newtheorem{theorem}{Theorem}

\begin{document}
\begin{frontmatter}
\title{\bf The tree-child network inference for line trees and the shortest common supersequences for permutations
}
\author[LB,LB2]{Laurent Bulteau}\ead{laurent.bultean@univ-eiffel.fr}
\author[LZ]{Louxin Zhang}\ead{matzlx@nus.edu.sg}

\address[LB]{Centre Nationale de Recherche Scientifique, FRANCE}
\address[LB2]{Université Gustave Eiffel, Paris, FRANCE}
\address[LZ]{Department of Mathematics, National University of Singapore,\\
10 Kent Ridge Road, Singapore 119076}


\begin{abstract} 
One strategy for inference of phylogenetic networks is to solve the  phylogenetic network problem, which involves inferring phylogenetic trees first and subsequently computing the smallest phylogenetic network that displays all the trees. This approach capitalizes on exceptional tools available for inferring phylogenetic trees from  biomolecular sequences.  Since the vast space of phylogenetic networks poses difficulties in obtaining comprehensive sampling, the researchers switch their attention to inferring tree-child networks from multiple phylogenetic trees, where in a tree-child network
each non-leaf node must have at least
one child that is an indegree-one node. 
  Two results are obtained in this work:  1) The tree-child network inference problem for multiple line trees remains NP-hard,  which is proved by a reduction from the shortest common supersequence problem for permutations. 2)  The tree-child networks with the least hybridization number that display all the line trees are the same as that display all the binary trees, whose hybridization number is  $\Theta(n^3)$  for $n (\geq 8)$ taxa.
\end{abstract}
\end{frontmatter}

\section{Introduction}

Recent genomic studies have highlighted the significant roles of recombination and introgression in genome evolution \cite{gogarten2005horizontal,koonin2001horizontal,Marcussen_14}. Consequently, there has been an increasing use of phylogenetic networks to model the evolution of genomes with  the presence of  recombination, introgression and other reticulate events \cite{Fontaine_15,koblmuller2007reticulate,Marcussen_14}. A phylogenetic network is a rooted directed acyclic graph (DAG) that represents taxa (genomes, individuals, or species) as its leaves and evolutionary events (speciation, recombination, or introgression) as its internal nodes. Over the past three decades, substantial progress has been made in understanding the theoretical aspects of phylogenetic networks \cite{gusfield2014book,huson2010book,steel2016phylogeny} (see also \cite{elworth2019advances,zhang2019clusters}).

The space of phylogenetic networks is vast, making it challenging to perform comprehensive sampling. As a result, popular methods like maximum likelihood and Bayesian inference, commonly used for phylogeny reconstruction, are not efficient enough for reconstructing phylogenetic networks containing a large number of reticulate events on more than 10 taxa \cite{lutteropp2022netrax,solis2016inferring,zhang2018bayesian}. This has prompted researchers to focus on inferring phylogenetic networks with specific combinatorial properties \cite{pickrell2012inference,van2022practical}. Popular classes of phylogenetic networks include galled trees \cite{gusfield2014book,wang2001perfect}, galled networks \cite{huson2009computing}, and tree-child networks \cite{cardona2009metrics2,cardona2020counting,zhang2019}. Furthermore, researchers are also investigating the parsimonious inference of phylogenetic networks from multiple trees, aiming to infer a network with the smallest hybridization number (HN) that display all the trees \cite{albrecht2012fast,mirzaei2015fast,wu2010close,yamada2020improved}, which we call the parsimonious networks. The HN, a generalization of the number of reticulate nodes in binary phylogenetic networks, quantifies the complexity of the network (refer to Section~\ref{sec2} for more details).

Inference of parsimonious phylogenetic networks is known to be NP-hard, even in the case of two input trees \cite{bordewich2007computing} and in the case tree-child networks are inferred \cite{linz2019attaching}.
 Notably,  a fast method has been recently developed to compute parsimonious tree-child networks for binary trees \cite{zhang2023fast}.

In this paper, using the approach developed in \cite{zhang2023fast} (summarized in Section 3), 
we study the inference of parsimonious tree-child networks for line trees.  We  prove that the inference problem remains NP-hard even for line trees (Sections 4-5). We also address the open problem of finding  the so-called ``universal" tree-child networks \cite{van2023three}. We show that the parsimonious tree-child networks for all line trees are identical to that display all binary trees (Section 6), for which a lower and upper bound for HN are given.

\section{Basic concepts and notation}
\label{sec2}

Let $X$ be a set of taxa.  A {\it phylogenetic network} on $X$  is a rooted DAG such that: 
\begin{itemize}
 \item The root is of indegree 0 and outdegree 1. There is at least one directed path from the root to every other node. \vspace{-0.5em}
    \item The leaves (which are of indegree 1 and outdegree 0) are  labeled one-to-one with the taxa.
    \vspace{-0.5em}
    \item All nodes except for the leaves and the root are either a {\it tree node} or a {\it reticulate node}. The former are of indegree 1 and outdegree 2, whereas the latter are of indegree more than 1 and outdegree 1. 
\end{itemize}
In a phylogenetic network, a node $u$ is said to be {\it below} another $v$ if there exists a directed path from $v$ to $u$.

A phylogenetic network is {\it binary} if every reticulate node is of indegree 2.
 A {\it binary phylogenetic tree} is a binary phylogenetic network that does not have any reticulate nodes. In this paper, a binary phylogenetic tree will be simply mentioned as a binary tree. A {\it line tree} is a binary tree in which all internal nodes but the  root  have  a leaf or two as  their children.

An important parameter of  phylogenetic networks is the {\it hybridization number} (HN).
For a phylogenetic network $N$,
${\rm HN}(N)=\sum_{v\in R(N)}(d_{in}(v) -1)$, where $R(N)$ is the set of reticulate nodes and $d_{in}(v)$ represents the indegree of $v$.
For a binary phylogenetic network $B$, each reticulate node has indegree 2 and thus ${\rm HN}(B)=\vert R(N)\vert$.

A {\it tree-child network} is a phylogenetic network in which  every non-leaf node has at least one child that is either a tree node or a leaf (Figure~\ref{def1}).

\subsection{The tree-child network problem }

Let $v$ be  a node of indegree 1 and outdegree 1 in a directed acyclic graph (DAG). Then, there is a unique edge $(u, v)$ entering 
$v$ and a unique edge $(v, w)$ leaving $v$ in the DAG. We may simplify it by removing $v$ and
replacing $(u, v)$ and $(v, w)$ with a new edge $(u, w)$. Such an operation is called the {\it degree-2 node contraction}.

A binary tree is {\it displayed} in a tree-child network if it can be obtained from
the network using the following two steps: (i) Delete all but one incoming edge for each reticulate
node. (ii) Contract  all the nodes with an indegree of 1 and an out-degree of 1.

We focus on 
 how to infer a tree-child network with  the minimum HN that displays all the input trees. This problem is formally defined as:
\begin{itemize}
\item[] {\bf The Tree-Child Network (TCN) Problem}
    \item[] {\bf Input} A set of binary trees on $X$. \vspace{-0.5em}
    \item[] {\bf Output} A  tree-child network with the minimum HN that displays all the trees.
\end{itemize}
The solution networks for the TCN problem are called {\it parsimonious}  tree-child network for the input trees.

\subsection{The shortest common supersequence  problem}

Let $A$ be an alphabet. A {\it string} on $A$ is an ordered sequence of characters. It is a {\it permutation string} if each character occurs exactly once in the string. 

A string on an alphabet is a {\it supersequence} of another if the latter can be obtained from the former by the deletion of 0 or more characters. A string is a {\it common supersequence} of multiple strings  if  it is a supersequence of every string. 

The length or size of a string $s$ is the total number of the occurrences of the characters in $s$, written as $\vert s\vert$.
A common supersequence is a shortest common supersequence (SCS) if it has the smallest length, over all the common supersequences  of the strings.
The SCS problem  is formally defined as:
\begin{itemize}
    \item[] {\bf Input} A set of strings on an alphabet. \vspace{-0.5em}
    \item[] {\bf Output} A SCS of the strings.
\end{itemize}
For a set $S$ of strings, every SCS string for the strings of $S$ has the same length. 
We will use ${\rm SCS}(S)$ to represent all the SCSs of the strings of $S$ and
$\vert{\rm SCS}(S)\vert$ to denote their length in the rest of this paper.
The SCS problem is a fundamental NP-complete problem \cite{garey1979computers}.

 \begin{figure}[th!]
\centering
\includegraphics[width=0.5\textwidth]{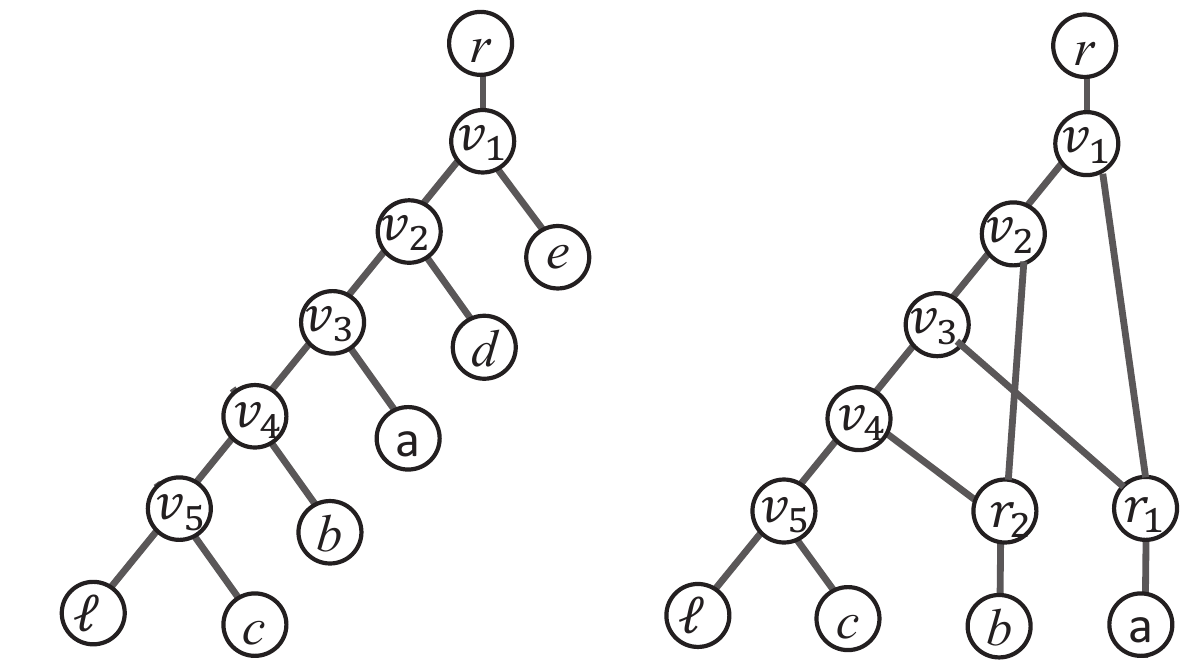}
\caption{\small 
The line tree defined using the permutation string $edabc$ (left) and 
the  tree-child network defined using by the string $ababc$ (right) on the set of characters appearing in the string and an extra symbol $\ell$. Here, the tree and network edges are oriented downwards and left-to-right.
The labels of non-leaf nodes are only used for the purpose of representing edges and they are not part of the tree and network.
\label{figure1}}
\end{figure}

\subsection{Permutation strings and one-component tree-child networks}

\begin{definition}
\label{def1}
Let $\Sigma$ be  an $n$-character alphabet and $\ell \notin \Sigma$.
For a permutation string $P=p_1p_2\cdots p_n$ on $\Sigma$,  $T(P, \ell)$ is defined to be the line tree on $\Sigma \cup \{\ell\}$ that has the node set $\Sigma \cup \{r, v_i, \ell\;\vert\; 1\leq i\leq n\}$ and the directed edge set $\{ (r, v_1), (v_i, v_{i+1}), (v_i, p_i), (v_n, \ell), (v_n, p_n) \;| \: 1\leq i\leq n-1\}$
{\rm (}left, Figure~\ref{figure1}{\rm )}.
\end{definition}
%


\begin{definition}
\label{def2}
Let $Q=q_1q_2\cdots q_m$ be a string on an $n$-character $\Sigma$ and 
$\ell\notin \Sigma$.  $N(Q, \ell)$ is defined to be the `one-component' tree-child network on $\Sigma \cup \{\ell\}$ that is obtained using the following two steps:

(i) Construct a DAG that has the node set
$\Sigma \cup \{r, \ell,  v_i, r_j\;|\; 1\leq i\leq m,  1\leq j\leq n\}$ 
and the directed edge  set $E_1\cup E_2$, where  $E_1=\{ (r, v_1), (v_i, v_{i+1}), (v_m, \ell), \vert 1\leq i\leq m-1\} \cup \{(r_j, a_j) \;|\; a_j\in \Sigma\}$ and $E_2$ contains $(v_i, r_j)$ if $q_i=a_j$ for every possible $i$ and $j$.

(ii) Contract all the nodes with an indegree of 1 and an outdegree of 1.
\end{definition}

Definition~\ref{def2} is illustrated in Figure~\ref{figure1} (right), where the node $r_3$ was removed in the node contraction step.
It is not hard to see that the HN of $N(Q, \ell)$ is $\vert Q\vert -n$, where $n=\vert \Sigma \vert$.

\section{Tree-child network inference via lineage taxa strings}

The parsimonious tree-child networks for multiple trees can be constructed from the lineage taxon strings (LTSs) of the taxa under an ordering on
$X$ \cite{zhang2023fast}. In this section, we shall restate the construction process on which our main results will be based.

Let $X$ consist of $n$ taxa and let $\pi: \pi_1<\pi_2 <\cdots <\pi_n$ be an ordering on 
$X$. 
We further assume that $\beta_1, \beta_2, \cdots, \beta_n$ are $n$  sequences satisfying the following conditions:
\begin{quote}
({\bf C1}) For each $i<n$, $\beta_i$ is a string on $\{\pi_{i+1}, \cdots, \pi_n\}$; 

({\bf C2}) $\beta_n$ is the empty sequence.
\end{quote}
It is proved in \cite{zhang2023fast} that the following algorithm outputs a tree-child network,
written as
$N\left(\pi, \{\beta_i\}^n_{i=1}\right)$, whose HN is equal to
$\sum_{1\leq i\leq n}\vert \beta_i \vert -n+1$. 
{\small 
\begin{center}
 \begin{tabular}{l}
 \hline
 {\sc Tree-Child Network Construction} \cite{zhang2023fast} \\
  1. ({\bf Vertical edges}) For each $\beta_i$,  define a 
  path $P_i$ with $\vert \beta_i\vert +2$ nodes: \\
  \hspace*{3em} $h_i, v_{i1}, v_{i2}, \cdots, v_{i\vert\beta_i\vert}, {\pi_i}$,\\
  \hspace*{1.5em}where $\beta_n$ is the empty sequence.\\
  2.  ({\bf Left--right edges}) 
  Arrange the $n$ paths
  from left to right as $P_1, P_2, \cdots, P_n$. \\ 
  \hspace*{1.5em}If the $m$-th symbol of $\beta_i$ is $\pi_j$, we add an
  edge $(v_{im}, h_{j})$ for each $i$ and  $m$. \\
  3. For each $i>1$, contract $h_i$ if $h_i$ is of indegree 1.   \\
    \hline
 \end{tabular}
 \end{center}
 } 
For example, applying the algorithm to $\beta_1=ababc$, $\beta_2=\beta_3=\beta_4=\epsilon$ and the ordering 
$\ell <a < b <c$,  we obtain the right tree-child network in 
Figure~\ref{figure1}.

Let $T$ be a binary tree on $X$.
For any $x, x'\in X$, we write 
$x<_{\pi} x'$ if $x$ is less than $x'$ under $\pi$. For a node $u$ of $T$, we use $\min_{\pi}(u)$ to denote the smallest of the taxa below $u$.  We label the root with the smallest taxon under $\pi$ and each non-root internal node $u$ with the larger of $\min_{\pi}(u')$ and $\min_{\pi}(u'')$, where $u'$ and $u''$ are the two children of $u$. 
In this way, the root and the remaining  $n-1$ internal nodes are uniquely labeled with a taxon. Moreover, the leaf $f$  is below the unique  internal node $w$ that had been labeled with $f$.  As a result, there exists a path $P_{wf}$ from $w$ to  $f$. The LTS  of the taxa $f$ consists of the taxon labels of the inner nodes in $P_{wf}$, ordered using the path orientation.

For example, if the alphabetic ordering
(i.e. $a<b<c<d<e<\ell$) is used, in the tree on left in Figure~\ref{figure1}, the root is labeled with $a$; $v_1$ to $v_5$ are labeled with $e, d, b, c, \ell$, respectively. Therefore, 
the LTS of $a, b, c$ are $edb, c, \ell$, respectively, whereas the LTS of $d, e, \ell$ are the empty string.

Consider $k$ binary trees $T_1, T_2, \cdots, T_k$ on $X$. 
 We write $\alpha_{ji}$ for the LTS of $\pi_i$ in  $T_j$ for each $i\leq n$ and each $j\leq k$. Then, for each $j$, $\alpha_{j1}, \alpha_{j2}, \cdots, \alpha_{jn}$ satisfy the conditions (C1) and (C2) (\cite{zhang2023fast}). Moreover, let $\gamma_i$ be a SCS of $\alpha_{1i}, \alpha_{2i}, \cdots, \alpha_{ki}$ for each $i$. 
 The sequences 
 $\gamma_1, \gamma_2, \cdots, \gamma_n$ also satisfy the conditions (C1) and (C2).

 \begin{theorem}
 \label{thm1}
    Let $T_j$ ($1\leq j\leq k$) be $k$ trees on an $n$-taxon $X$ and let $P$ be a tree-child network on $X$ that displays all the trees.
    If $P$ is a parsimonious tree-child network for $\{T_1, T_2, \cdots, T_k\}$,  then   an ordering  $\pi$ on $X$ can be computed in linear time such that:

    {\rm (i)}
    $P=N\left(\pi,
    \{\gamma_i
    \}^n_{i=1}\right)$, where  $\gamma_i$ is a SCS of $\alpha_{1i}, \alpha_{2i}, \cdots, \alpha_{ki}$, where  $\alpha_{ji}$ is the LTS of $\pi_i$  under $\pi$ in the tree $T_j$, 
    
    {\rm (ii)} The SCS strings 
    $\gamma_1, \gamma_2, \cdots, \gamma_n$ can be computed by labelling the internal nodes of $P$ in  linear time, 
    and
    
    {\rm (iii)}  ${\rm HN}(P)=\sum_{1\leq i\leq n}\vert \gamma_i\vert -n+1
    = \sum_{1\leq i\leq n}\vert {\rm SCS}(\alpha_{1i}, \alpha_{2i}, \cdots, \alpha_{ki})\vert -n+1$.
 \end{theorem}
 
The proof of Theorem~\ref{thm1} appears in Section A of the Supplemental Methods of \cite{zhang2023fast}. 
Theorem~\ref{thm1} implies that the TCN problem has usually multiple solution networks with the same HN.

\section{Equivalence of the TCN and SCS problems}
 
According to Theorem~\ref{thm1},  the TCN problem can be solved by reducing it to multiple SCS sub-problems with instances being the LTSs of  taxa through examining all possible orderings on the taxa. 
To establish a reduction from the SCS problem to the TCN problem, we show that any given SCS instance—comprising a collection $P$ of permutation strings on $\Sigma$—can be efficiently transformed into a corresponding TCN instance $I_P$ such that
for $s\in {\rm SCS}(P)$,  $N(s, \ell)$ is a solution to $I_P$, $\ell \notin \Sigma$.

Consider an instance of the SCS problem with the input set $P$ consisting of $k$ permutation strings $P_i$ ($1\leq i\leq k$) on an $n$-character alphabet $\Sigma$. 
%
By Theorem~\ref{thm1}, each parsimonious tree-child network $N$ for the $k$ trees $T(P_1, \ell), T(P_2, \ell), \cdots, T(P_k, \ell)$ can be constructed from the LTSs of taxa under some ordering  $\pi: \pi_1< \pi_2 <\cdots <\pi_n <\pi_{n+1}$ on $\Sigma\cup\{\ell\}$, where $\ell \not\in \Sigma$.  Importantly, $\pi$ can be found in linear time according to the theorem. We now prove that  a SCS for $P$ can be obtained  from the SCSs of the LTSs of the taxa under $\pi$.  The latter are found in $N$ (Theorem~\ref{thm1}). 
In the rest of the discussion, 
we consider the following two cases depending on whether $\ell$ is the smallest  taxon under $\pi$ or not. 

{\bf Case 1}. $\ell=\pi_1$. 

By definition, in each $T(P_i, \ell)$,  the LTS of $\ell$ is $P_i$  and empty for every other  $\pi_i$,  $i>1$. In this case, 
by Theorem~\ref{thm1}.ii, a SCS $\gamma_1$ of $P_1, P_2, \cdots, P_k$ is computed from the parsimonious network $N$ in linear time.

{\bf Case 2}. $\ell=\pi_t$, where $t>1$. 



For each $i$, we let $P_i=p_{i1}p_{i2}\cdots p_{in}$, where $p_{ij}\in \Sigma$ for each $j$.
We first identify which taxa have a non-empty LTS under $\pi$ in the line tree $T(P_i, \ell)$.
They can be found as follows.

Initially, define $\beta_{i1}=\pi_1$. Assuming that we have obtained $\beta_{ij}=p_{ix} <_{\pi} \min \{p_{in}, \ell\}$, 
we compute $\beta_{i (j+1)}=\min_{\pi} \{ p_{i(x+1)}, \cdots, p_{in}, \ell \}$.
Keep on this procedure, we obtain a finite sequence of taxa: 
\begin{eqnarray}
\label{eqn1_LTS}
    \beta_{i1}=\pi_1, \beta_{i2}, \cdots, \beta_{iw_i}=
\min _{\pi} \{p_{in}, \ell\},
\end{eqnarray}
such that 
$\beta_{ik}<_{\pi} \beta_{i(k+1)}$  and the leaf 
$\beta_{i(k+1)}$  is deeper than $\beta_{ik}$ 
in $T(P_i, \ell)$.

\begin{figure}[t!]
\centering
\includegraphics[width=0.8\textwidth]{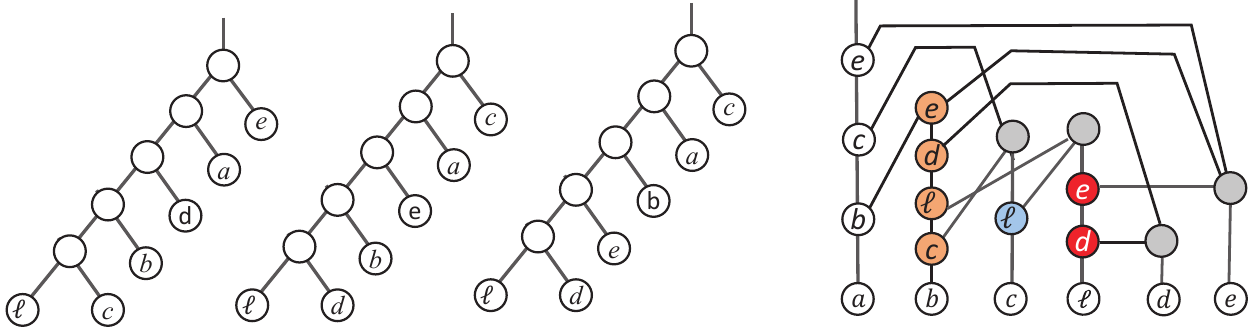}
\caption{\small 
Illustration of the proof in Case 2.
The three trees $T_1, T2$ and $T_3$ (left to right) are defined using
the permutations $P_1=eabdc, P_2=caebd, P_3=cabed$, respectively. The tree-child network on right is constructed from the three trees using {\sc Tree-Child Network Reconstruction} under the  ordering $a<b<c<\ell<d<e$. In the trees and  network,  edges are oriented downward and from left to right and the root is not drawn. The colored paths corresponds with the nonempty  LTSs of 
$a, b, c, \ell$, respectively.
\label{figure2}}
\end{figure}

For example, we consider the trees in Figure~\ref{figure2}. Under the ordering:
$a<b<c<\ell<d<e$, in the left tree $T_1$, we have:
$$\beta_{11}=a, \;\;\beta_{12}=b, \;\; \beta_{13}=c=\min \{c, \ell\}, \;\; w_1=3.$$
On the other hand, under the same ordering, in the middle tree $T_2$, we have:
$$\beta_{21}=a, \;\;\beta_{22}=b, \;\; \beta_{23}=\ell=\min \{c, \ell\}, \;\; w_1=3.$$

By the definition of LTS,  the LTS of  $\beta_{ij}$ ends with $\beta_{i(j+1)}$ and thus are nonempty under $\pi$ for each $j< w_i$.  Note that $p_{in}$ and $\ell$ are siblings in the bottom of the tree. If $\beta_{w_i}=p_{in}$, 
the LTS of $\beta_{iw_i}$ ends with  
 $\ell$. Otherwise, $\beta_{w_i}=\ell$ and the LTS of $\beta_{iw_i}$ ends with $p_{in}$  under $\pi$.
 In addition, the LTS is empty  for any other taxon under $\pi$.  
 
 Continue the discussion on $T_1$ in Figure~\ref{figure2}, 
 The LTS of $a$, $b$ and $c$ are 
 $eb$, $dc$, and $\ell$, respectively. The LTS of $e, d, \ell$ are empty.

 Furthermore,   we have the following fact, whose proof is straightforward.

\begin{proposition}
\label{prop1}
    Let the LTS of $\beta_{ij}$ defined in Eqn.~(\ref{eqn1_LTS}) be
    $S_{ij}$ in $T_i$ under $\pi$, where $ \pi_1 <_{\pi} \ell$.  Then, for each $i$,
    \begin{eqnarray*}
    && P_i=
     S'_{i1}\beta_{i1}S'_{i2}\beta_{i2}\cdots 
    S'_{i(w_i-1)} \beta_{i(w_i-1)}S''_{iw_i}, \\
    &&S''_{iw_i}=\left\{\begin{array}{ll}
    S_{iw_i} & \mbox{ if } \beta_{iw_i}=\ell \\
    S'_{iw_i} \beta_{iw_i} & \mbox{ if } \beta_{iw_i} \neq \ell \end{array}\right.
    \end{eqnarray*}
    where 
    $S'_{it}$
    denotes the string obtained from $S_{it}$ by removal of its last character for any $t$
    and the right-hand side is the concatenation of the strings and characters. 
\end{proposition}


Let us the notation in Prop.~\ref{prop1} for further discussion. 
Since, by assumption, $N$ is a parsimonious tree-child network  for $T(P_1, \ell), T(P_2, \ell), \cdots, T(P_k, \ell)$,
and $\pi$ is an ordering satisfying the conditions 
(i)--(iii) in Theorem~\ref{thm1}. 
We further assume 
$\gamma_i$ is a SCS of the LTSs of $\pi_i$ found in $N$ for each $i$ (Theorem~\ref{thm1}.(ii)).  From the above discussion, we have the following fact.

\begin{proposition}
 Each SCS string $\gamma_j$ is non-empty if and only if 
 $\pi_j=\beta_{it}$,
 for some $i$ and $t\leq w_i$ appearing in Eqn.~(\ref{eqn1_LTS}).
\end{proposition}

For example, for the network in Figure~\ref{figure2}, the nonempty SCSs are $\gamma_a$ (white), $\gamma_b$ (orange), $\gamma_c$ (blue) and $\gamma_\ell$ (red).

For each $i$ such that $\gamma_i$ is nonempty, we define $\gamma'_i$ to be the string obtained $\gamma_i$ by removal of its last character.
Assume that $\gamma'_{h_1}, \gamma'_{h_2}, \cdots, \gamma'_{h_j}$ be all the obtained strings from nonempty $\gamma_i$'s, where $$\pi_{h_1}<_\pi 
   \pi_{h_2}<_\pi \cdots <_{\pi} \pi_{h_j}.$$
We further define 
$Q$ to be the string 
obtained from:
\[w_{h_1}\pi_{h_1}w_{h_2}\pi_{h_2}\cdots w_{h_{j-1}}\pi_{h_{j-1}} w_{h_j}\pi_{h_j}\] by removal of all the occurrences of $\ell$ if any, where, for each $t\leq j$, 
$$w_{h_t}=\left\{\begin{array}{ll}
 \gamma'_{h_t} & \mbox{ if } h_t\neq \ell,\\
 \gamma_{h_t} & \mbox{ if } h_t=\ell.
\end{array}\right.$$

\begin{figure}[t!]
\centering
\includegraphics[width=0.6\textwidth]{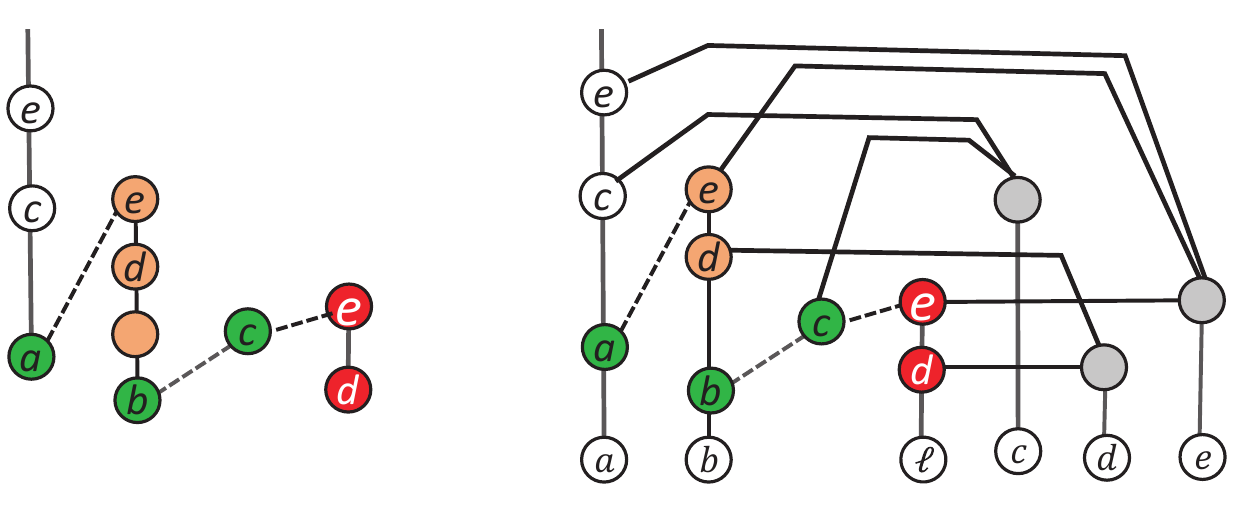}
\caption{\small  The string $Q=ecaedbced$ (left) and the  tree-child network $N(Q, \ell)$ (right), which are derived from the network in Example 1. The dashed lines are used to connect the substrings/paths obtained from different LTSs. The non-labeled orange node had $\ell$ as its label in the original network and was removed to compute $Q$.   The network edges are drawn downwards and from left to right and the root is not drawn. 
\label{figure3}}
\end{figure}

For the example in Figure~\ref{figure2}, the computation of the string $Q$ is illustrated in Example 1 (below) and in Figure~\ref{figure3}.

Since $\gamma_{h_j}$ is a SCS of the LTSs of $\pi_{h_j}$ in each tree, by Proposition~\ref{prop1}, 
$Q$ is a super-sequence of $P_1, P_2, \cdots, P_k$ and 
   \begin{equation}
       \label{eqn2_lts}
   {\rm HN}(N(Q, \ell)) =
   \vert Q\vert -n \leq 
   \sum _{1\leq j\leq n} |\gamma_i| -n 
   ={\rm HN}(N).
   \end{equation}

  Since $N$ is a parsimonious tree-child network for the line trees, 
  \begin{equation}
  \label{eqn3_lts}
  {\rm HN}(N)\leq {\rm HN}(N(\gamma, \ell))=\vert {\rm SCS}(P)\vert -n,
  \end{equation}
  where $\gamma$ is a SCS of the permutation strings of $P_1, P_2, \cdots, P_k$.
  

  Taken together, Inequalities~(\ref{eqn2_lts}) and
  (\ref{eqn3_lts}) implies that $Q$ constructed above from $N$ is a SCS of the permutation strings 
  $P_1, P_2, \cdots, P_k$.
  This proves the following result.

  \begin{theorem}
      \label{thm2}
    Let $N$ be the parsimonious tree-child network for $T(P_1, \ell), T(P_2, \ell)$, $\cdots$, $T(P_k, \ell)$. We can conmput a SCS of the strings of $P$ in linear time form $N$.
    \end{theorem}

  \noindent {\bf Example 1}. Consider the ordering $\pi: a<b<c<\ell<d<e$ for the tree lines trees in Figure~\ref{figure2}.  The LTSs of the taxa under $\pi$ in the three trees are listed in the following table.
   \\

   \begin{tabular}{c|ccc|c}
   \hline
     Taxon   &  LTS in $T(P_1)$ &  LTS in $T(P_2)$ &  LTS in $T(P_3)$ & SCS\\
     \hline
      $a$  & $eb$ & $cb$ & $cb$ & $ecb$\\
      $b$ &  $dc$ & $e\ell$ & $\ell$ & $ed\ell c$\\
      $c$ &  $\ell$ & $\epsilon$ & $\epsilon$ &$\ell$\\
      $\ell$ & $\epsilon$  & $d$ & $ed$ & $ed$\\
      $d$ & $\epsilon$ &$\epsilon$ &$\epsilon$ & $\epsilon$\\
      $e$ &$\epsilon$ &$\epsilon$ & $\epsilon$& $\epsilon$\\
      \hline
   \end{tabular}
   \\
   
  \noindent In the table, $\epsilon$ denotes the empty string. 
  The last column lists a SCS of the LTSs $\gamma_i$ for each taxon under $\pi$. 
  Since nonempty SCSs are $\gamma_a, \gamma_b, \gamma_c, \gamma_\ell$, 
 the string $Q$ constructed before Inequality~(\ref{eqn2_lts}) is $ecaedbced,$
   which is obtained from $\gamma'_aa\cdot \gamma'_bb\cdot \gamma'_cc\cdot \gamma'_{\ell}\ell=(ec)a (ed\ell) b c (ed)\ell$  after removing two $\ell$'s.
Clearly, $Q$ is a supersequence of $eadbc, caebd$, and $cabed$. 
The  tree-child network $N(Q)$ is shown  in Figure~\ref{figure3} and has a HN of 4. On contrast, the HN of the network in Figure~\ref{figure2} is 5, implying that the latter is not a parsimonious. 

\section{NP-hardness of the SCS problem for permutation strings}

  SCS is already known to be NP-hard when all input strings consist of 2 distinct characters~\cite{timkovskii1989complexity}.
  Let us denote this variant {\em 2-SCS} (we further need the trivial constraint that no character appears in every input string).
  We use this fact to obtain the following NP-hardness result.
  
\begin{theorem}
\label{thm3}
    The SCS problem is NP-hard even for permutation strings.
\end{theorem}
\noindent {\bf Proof.}
   Consider an instance $\mathcal S$ of 2-SCS with $m$ length-2 strings over a size-$n$ alphabet  $X=\{x_1,\ldots,x_n\}$, and an integer $k$.  Let $N=n+k+1$, and create a size-$N$ set $Y=\{y_1,\ldots,y_N\}$ of \emph{separators}. In the context of strings, we also write $X$ and $Y$ for the strings $x_1\ldots x_n$ and $y_1\ldots y_N$, respectively. For any string $ab\in \mathcal{S}$ (with $a,b\in X$ and  $a\neq b$), we write $X_{-ab}$ for the subsequence of $x_1x_2\ldots x_n$ obtained by removing $a$ and $b$, and $S_{ab} = a b \cdot Y \cdot X_{-ab}$. Note that each $S_{ab}$ is a permutation string on $X\cup Y$.  Let us write $\mathcal S' = \{S_{ab}, ab\in\mathcal S\}$ and $k'=k+N+n$. We now prove the following equivalence that completes the reduction.

    \begin{center}
        Strings in $\mathcal S$ have a common supersequence $T$ of size $k$ \\$\Leftrightarrow$ Strings in $\mathcal S'$ have a common supersequence $T'$ of size $k'$
    \end{center}

    $\Rightarrow$ 
      Build $T' = T \cdot Y \cdot X$.  $T'$ is a length-$k'$ string, and it is a supersequence of any $S_{ab}$ for $ab\in S'$ (since $T$ is a supersequence of $ab$ and $X$ is  a supersequence of $X_{-ab}$).        

    $\Leftarrow$
      Pick such a string $T'$. It contains at least one occurrence of $Y$ as a subsequence. Let $P,R$ be the matching prefix and suffix of $T'$ (i.e. $T'=P\cdot R$) such that $R$ is the smallest suffix containing $Y$ as a subsequence.
      Let $T$ be the subsequence of $P$ obtained by removing all separator characters.  
      We have $|P|\leq k'-N = k+n <N$, so $P$ may not contain an entire copy of $Y$. Hence, for any $S_{ab} = ab\cdot Y\cdot X_{-ab}\in \mathcal S'$, we have that $ab$ is a subsequence of $P$ and $X_{-ab}$ is a subsequence of $R$.
      Overall, $P$, and  also $T$, are common supersequence of all $ab\in \mathcal S$, and $R$ is a common supersequence of all $X_{-ab}$. In order to bound their sizes, note that $R$ contains each character of $X$ and $Y$ at least once, so $|R|\geq N+n$. 
      Hence, $T$ has size at most $k'-N-n=k$, and is a common supersequence of $\mathcal S$.
      This concludes the proof.
\\

Taken together, Theorem~\ref{thm2} and Theorem~\ref{thm3} imply the following theorem.

\begin{theorem}
    The TCN problem is NP-hard even for line trees. 
\end{theorem}

\noindent {\bf Open problem 1.} Does the TCN problem remain NP-hard for two line trees?

The TCN problem for three line trees was studied by Van Iersel et al. in \cite{van2023three}.

\section{Tree-child networks that display all the line trees}

\begin{theorem}
\label{thm5}
  The parsimonious tree-child networks that display all the line trees  on the taxon set $X$ are identical to that that display all the binary trees. 
\end{theorem}
{\bf Proof.} 
Let $\pi: \pi_1 < \pi_2<\cdots < \pi_n$ be an arbitary ordering on $X$. Set 
$X_i=\{\pi_i, \pi_{i+1}, \cdots, \pi_n\}$.
Clearly, $X=X_1$. We also use $X_i^k$ to denote the set of all the strings of length $k$ in which each character appears no more than once, where $k\geq 1$.
Let $S_i$ be the set of the LTSs of $\pi_i$ obtained in all the line trees.
For each length-($n-1$) string $s=a_1a_2\cdots a_{n-1}$ on $X_2$, 
we obtain the line $T(s, \pi_1)$.  Since $\pi_1$ is one of the deepest leaf, 
we conclude that $s$ is the LTS of $\pi$ 
in $T(s, \pi_1)$. This  implies that 
$X_2^{n-1}\subseteq S_1$. Since the LTS of $\pi_1$ in any line tree is a string in which each character of  $X_2$ appears at most once, 
\begin{equation}
\label{eqn4}
    S_1\subseteq X_2 \cup X_2^2 \cup \cdots 
\cup X_2^{n-1}
\end{equation}

Similarly, we can prove that, for each $k\leq n-1$, 
\begin{equation}
\label{eqn5}
X_{k+1}^{n-k} \subseteq S_k\subseteq X_{k+1} \cup X_{k+1}^2 \cup \cdots 
\cup X_{k+1}^{n-k}.
\end{equation}
Taken together, Eqn.~\ref{eqn4} and
\ref{eqn5} implies that for each $k$,
\begin{equation}
    \vert {\rm SCS}(S_k)\vert =\vert {\rm SCS}(X_{k+1}^{n-k}) \vert
\end{equation}

In the case for all the trees, Eqn. (\ref{eqn4}) remains for all $k$. Similarly,  Eq. (\ref{eqn5}
) is also true. This implies that for any ordering on $A$, the LTSs of any specific taxon in all the lines and in all binary trees have the same SCS. Therefore, the parsimonious trees for all line trees are identical to that for all binary trees.
This concludes the proof.
\\

\begin{lemma}
\label{lemma1}
    Let $A$ be a size-$n$ alphabet and let $A_p$ denote the set of all permutation strings on $A$. Then, 
     \begin{eqnarray*}
   && \vert \mbox{\rm SCS}(A_p) \vert = \left\{ \begin{array}{cc}
     1,  &  \mbox{if } n=1; \\
     3,  &  \mbox{if } n=2;\\
     n^2 -2n+4,&  \mbox{ if } 3\leq n\leq  7;
   \end{array} \right.\\
   &&   n(n+1)/2 \leq \vert \mbox{\rm SCS}(A_p)\vert \leq 
    n^2-2n+4- \lfloor (n-7)/3\rfloor,  \mbox{ if } n \geq 8.
     \end{eqnarray*} 
\end{lemma}
{\bf Proof.}
The results for $n\leq 7$ can be verified using computing. For example, the string $1213121$ is a SCS of all the permutation strings on $\{1, 2, 3\}$. The string 
$123413214321$ a SCS of all the permutation strings on $\{1, 2, 3, 4\}$. 
The length-19 string $1234512341523142351$ is a SCS of all the permutation strings on $\{1, 2, 3, 4, 5\}$. 

For $n\geq 8$, the lower bound was given by  Kleitman and Kwiatkowski \cite{kleitman1976lower}  and the upper bound was proved by  Radomirovi\'c \cite{radomirovic2012construction}.
\\

\begin{theorem}
\label{thm5}
Let $P$ be a parsimonious tree-child network that display all the binary trees  on the taxon set $X$. Then, 
  $(n-2)(n-1)(n+3)/6\leq {\rm HN}(P) \leq (n^3-5n^2+21n-53)/3$, where $n=|X|\geq 8$.
\end{theorem}
{\bf Proof.} 
Let $f(x)=\lfloor (x-7)/3\rfloor$.
By Lemma~\ref{lemma1},  $\vert {\rm SCS}(S_k)\vert =\vert {\rm SCS}(X_{k+1}^{n-k}) \vert \leq 
(n-k)^2-2(n-k)+4-f(n-k)$ for each $1\leq k\leq n-8$,  $\vert {\rm SCS}(S_k)\vert=(n-k)^2-2(n-k)+4$ for $n-7\leq k \leq n-3$, $\vert {\rm SCS}(S_{n-2})\vert =3$ and $\vert {\rm SCS}(S_{n-1})\vert =1$. Therefore, 
the HN of the tree-child network $Q_\pi$ constructed using {\sc Tree-Child Network Construction} from ${\rm SCS}(S_k)$ ($1\leq k\leq n-1$)  is:
\begin{eqnarray*}
  {\rm HN}(Q_{\pi}) &\leq & 
    4+ \sum^{n-3}_{k=n-7} [(n-k)^2-2(n-k)+4]\\
    & & 
     + \sum^{n-8}_{k=1} [(n-k)^2-2(n-k)+4-f(n-k)] -(n-1)\\
     &=& 5-n +\sum^{n-1}_{j=3}(j^2-2j+4) 
     -\sum^{n-1}_{j=8}f(j)\\
     &=& 3(n-2)+ n(n-1)(2n-7)/6 -\sum^{n-1}_{j=8}\lfloor (j-7)/3\rfloor\\
\end{eqnarray*}
If $n=3k+2, 3k$, $k\geq 3$, the last sum is $(n-8)(n-9)/6$. If $n=3k+1, k\geq 3$, the last sum term is 
$(n-7)(n-10)/6$. 
Therefore, for $n\geq 8$, 
$$ {\rm HN}(Q_{\pi})
\leq (n^3-5n^2+21n-53)/3. $$

According to Lemma~\ref{lemma1},  $\vert {\rm SCS}(S_k)\vert =\vert {\rm SCS}(X_{k+1}^{n-k}) \vert \geq (n-k)(n-k+1)/2$ for $1\leq k\leq n-1$. 
The lower bound on $ {\rm HN}(Q_{\pi})$ is obtained as:
\begin{eqnarray*}
    {\rm HN}(Q_{\pi}) & =&
    \sum^{n-1}_{k=1} \vert {\rm SCS}(S_k)\vert  -(n-1)\\
    &\geq &  \sum^{n-1}_{k=1} (n-k)(n-k+1)/2 -(n-1)\\
    &=&  (n-2)(n-1)(n+3)/6.
\end{eqnarray*}
This concludes the proof.
\\

\begin{figure}[t!]
\centering
\includegraphics[width=0.3\textwidth]{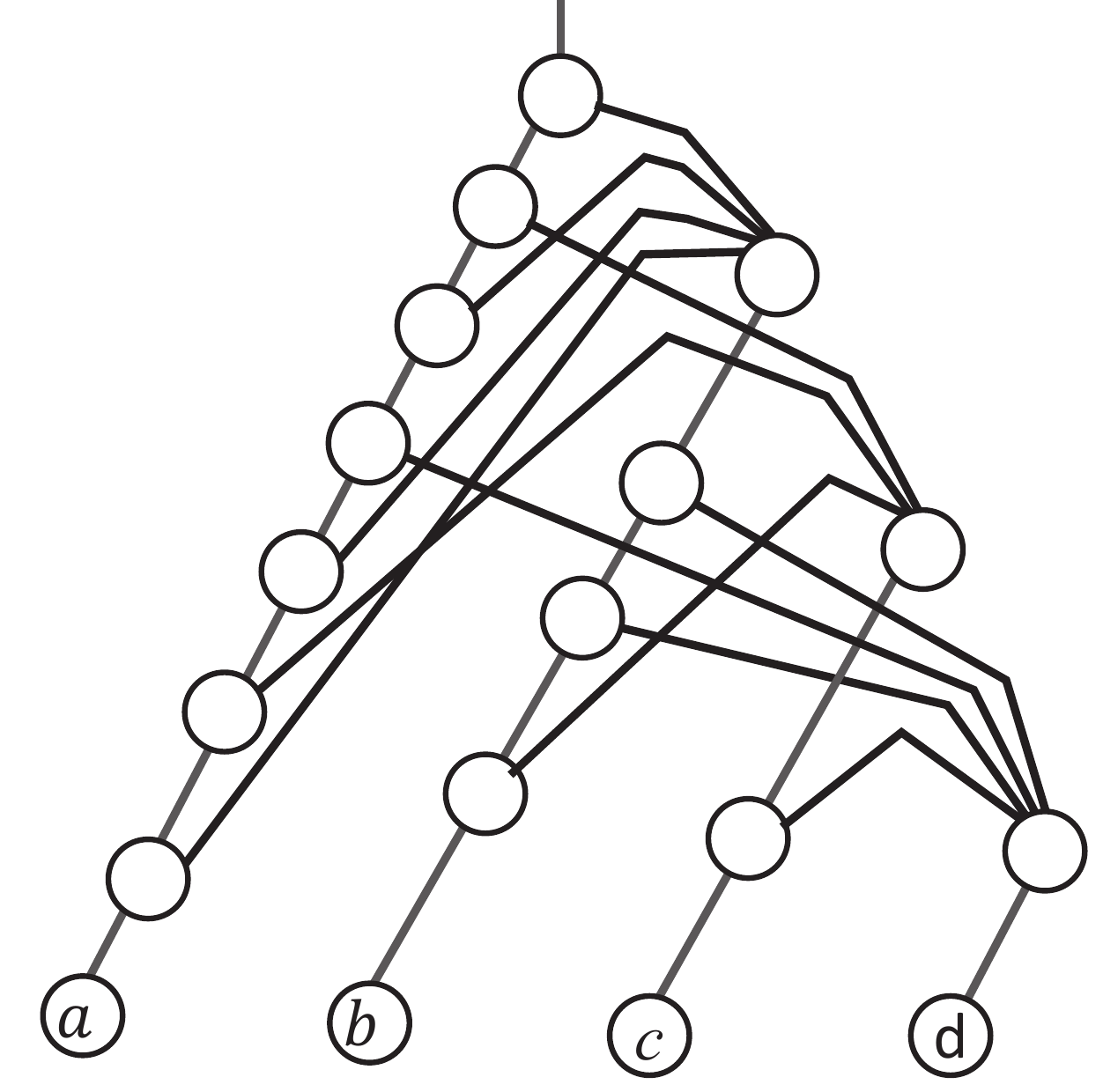}
\caption{\small  A parsimonious tree-child network with $8$ reticulate nodes that displays all 15 binary trees on taxa $a, b, c, d$. The straight edges are oriented downwards and other edges are oriented from left to right.
\label{fig4}}
\end{figure}

\noindent {\bf Example 2}. For $n\leq 5$,  a parsimonious tree-child network for all the 15 binary trees on $\{a, b, c, d\}$ is given in Figure~\ref{fig4}. 

\section{Conclusions}

In this work, we have made two contributions to the study of tree-child networks. Firstly, the tree-child network inference problem remains NP-hard, 
further fortifying the NP-hardness result originally presented in \cite{linz2019attaching}. 
To prove the NP-hardness, we prove that the shortest common supersequence problem is also NP-hard, which holds intrinsic interest for the broader scientific community.

Secondly, we have proved that the HN of the so-called "universal tree-child network" on $n$ taxa is $\Theta(n^3)$, illuminating the expressive power of tree-child networks. 

\section*{Acknowledgements} 
LX Zhang was partially supported by  Singapore
MOE Tier 1 grant R-146-000-318-114 and Merlin 2023. He thanks Yufeng Wu for useful discussion in the early stage of this work.


\end{document}